\documentclass[11pt]{amsart}
\usepackage{amsmath, amssymb}
\usepackage{amsthm}
\usepackage{listings}
\usepackage{tikz-cd}
\newtheorem{thm}{Theorem}[section]
\newtheorem{lem}[thm]{Lemma}

\newtheorem{prop}[thm]{Proposition}

\newtheorem{rem}[thm]{Remark}
\newtheorem*{thm*}{Theorem \ref{th:main3}}
\newtheorem*{thm*2}{Theorem \ref{th:main1}}
\usepackage[margin=1.66in]{geometry}

\theoremstyle{definition}

\renewcommand{\rho}{\varrho}

\newcommand{\Sig}{\mathbb{S}}

\newcommand{\pr}[1]{\mathtt{#1}}

\date{}

\begin{document}

\title{Definability of some $k$-ary Relations Over Second Order kinds of Logics}
\author[S. Costa]{Simone Costa}
\address{DICATAM - Sez. Matematica, Universit\`a degli Studi di Brescia, Via Branze 43, I-25123 Brescia, Italy}
\email{simone.costa@unibs.it}
\author[M. Dalai]{Marco Dalai}
\address{DII, Universit\`a degli Studi di Brescia, Via Branze 38, I-25123 Brescia, Italy}
\email{marco.dalai@unibs.it}
\author[S. Della Fiore]{Stefano Della Fiore}
\address{DII, Universit\`a degli Studi di Brescia, Via Branze 38, I-25123 Brescia, Italy}
\email{stefano.dellafiore@unibs.it}
\author[A. Pasotti]{Anita Pasotti}
\address{DICATAM - Sez. Matematica, Universit\`a degli Studi di Brescia, Via Branze 43, I-25123 Brescia, Italy}
\email{anita.pasotti@unibs.it}

\begin{abstract}

We consider the exprissibility  in monadic second order logic of certain relations of importance in computer science. For integers $n\geq 1$ and  $k\leq b$, a $k$-tuple of sequences in $\{0,1,\ldots, b-1\}^n$ are said to be $k$-hashed if there is a coordinate where they all differ. A set $\mathcal{C}$ of sequences is said to be a $k$-hash code if any $k$ distinct elements are $k$-hashed. Testing whether a code is $k$-hashing and determining the largest size of $k$-hash codes is an important problem in computer science. The use of general purpose solvers for this problem leads to question what minimal logic is needed to represent the problem.

In this paper, we prove that the $k$-hashing relation on $k$-tuples is not definable in Monadic Second Order Logic (MSO), highlighting its limitations for this problem. Instead, the property can be expressed in extensions of the MSO that add the equi-cardinality relation. 

\end{abstract}

\maketitle

\section{Introduction}

Let $\mathcal{C}$ be a subset of $\{0,1,\ldots,b-1\}^n$ with the property that, for any $k\leq b$ distinct elements, there exists a coordinate in which they all differ. A code $\mathcal{C}$ with this property is called a $(b,k)$-hash code of length $n$, and plays an important role in computer science. In particular,
as $n$ grows, of relevant theoretical interest is the determination of the cardinality $T_{b,k}(n)$ of the largest such $(b,k)$-hash codes \cite{KoS, KM}. Upper and lower bounds were derived over the years, showing that the particular case $b=k=3$, while being the simplest non-trivial instance to formulate, actually represents the most challenging one. A simple polynomial upper bound of the form $T_{3,3}(n)\leq C (3/2)^n$, for some constant $C$, had remained for years the best known (see \cite{CD}, \cite{DFGP} or \cite{KS} for discussions), and only recently a polynomial improvement of the form $T_{3,3}(n)\leq C n^{-2/5} (3/2)^n$ was derived in \cite{B}. On the achievabiliy side, the best known lower bound $T_{3,3}\geq C (9/5)^{n/4}$ can be derived by means of an explicit constructions for $n=4$ (\cite{KM}; see also \cite{bis-dha-gij-2024} for an alternative method).

In this context, testing for arbitrary $n$ and $T$ whether there exist $(b,k)$-hash codes of length $n$ and size $T$ represents a relevant question, and explicit  algorithms for answering it are of interest (\cite{DFGP}, \cite{KS}). One naturally derived question is whether general purpose solvers can be used for this task. Methods for proving the existence (or the non-existence) of structures that satisfy given constraints have been recently implemented, as for example in the Satisfiability Modulo Theories (SMT, for short) program MONA \cite{Mona}, where these constraints can be written in the First Order Logic or in the Monadic Second Order Logic (MSO).

Our attempts to use these approaches to improve the results of \cite{DFGP} and \cite{KS} lead us to ask ourselves whether the $(b,k)$-hashing $k$-ary relation is even expressible in the Monadic Second Order Logic. In this paper, we show that the answer is negative.

In the interest of exposition, we consider explicitly the case $b=k=3$, the general case being easily derivable from this one. In this special case, three sequences having the $(3,3)$-hashing property are simply said to be \emph{trifferent}, and we thus present our method by showing that the \emph{trifference} relation for general length $n$ cannot be expressed in MSO logic. In particular, in Section \ref{sec:3}, we will show that there are no formulas $f(X,Y,Z)$ of rank $\rho$ satisfied if and only if $X,Y,Z$ are trifferent words of length large enough with respect to $\rho$.  Firstly, in Subsection \ref{subsec:33}, we will provide a  proof of the fact that there are no formulas $f(X,Y,Z)$ satisfied whenever $X,Y,Z$ are finite, trifferent words (i.e. the trifference relation is not S3S definable). Then we prove the result more in general for both finite and infinite words.
On the other hand, if we extend the MSO by adding the possibility of expressing the equality between cardinalities, it is possible to characterize the trifference property between three words. There are, indeed, several well-studied theories on this topic such as the counting MSO (CMSO, in brief), \cite{C1, C2}, the MSO with cardinalities \cite{Felix}, and extensions which add a weak kind of arithmetics (the so-called Presburger arithmetics), see \cite{K1,K2, Rudolph}.
In the last section of this work, we will provide a program in the SMT solver Z3 \cite{Z3} that can be applied to show the non-existence of a trifferent code of size $11$ and length $5$ (which recovers a result of \cite{DFGP}).

\section{Preliminars}
Monadic Second-Order logic is a widely used and expressive logical framework that balances expressive power with computational tractability. It is particularly well-suited for analyzing finite or countable structures. MSO extends beyond the capabilities of first-order logic by enabling the expression of  ``mildly recursive'' structural properties, such as connectedness and reachability, which are essential for addressing key modeling needs in areas like database theory and knowledge representation.

The \textbf{monadic second-order theory} $S3S$ is defined over a specific class of structures representing infinite ternary trees. Formally, the structures considered in $S3S$ are of the form:
$$
\mathfrak{T} = (T, \pr{r},\succ_0^{\mathfrak{T}}, \succ_1^{\mathfrak{T}}, \succ_2^{\mathfrak{T}}),
$$
where:
\begin{itemize}
\item $T = \{0,1,2\}^*$ is the set of all finite ternary strings, serving as the \emph{domain} of the structure, and representing the nodes of a complete infinite ternary tree;
\item $\pr{r}$ is the root of the tree, i.e. the node associated to the empty string;
\item $\succ_0^{\mathfrak{T}}$ is a subset of $T \times T$ called the \emph{left successor relation}, defined as:
$$
\succ_0^{\mathfrak{T}} = \{(w, w0) \mid w \in \{0,1,2\}^*\};
$$
\item $\succ_1^{\mathfrak{T}}$ is a subset of $T \times T$ called the \emph{middle successor relation}, defined as:
$$
\succ_1^{\mathfrak{T}} = \{(w, w1) \mid w \in \{0,1,2\}^*\};
$$
\item $\succ_2^{\mathfrak{T}}$ is a subset of $T \times T$ called the \emph{right successor relation}, defined as:
$$
\succ_2^{\mathfrak{T}} = \{(w, w2) \mid w \in \{0,1,2\}^*\}.
$$
\end{itemize}
It is straightforward to define, in analogy, the theory $SbS$ over the infinite $b$-ary tree. Formally, the structures of $SbS$ are of the form
\[
\mathfrak{T} = \bigl(T, \pr{r}, \succ_0^{\mathfrak{T}}, \succ_1^{\mathfrak{T}}, \dots, \succ_{b-1}^{\mathfrak{T}} \bigr),
\]
where:
\begin{itemize}
    \item $T = \{0,1,\dots, b-1\}^*$ is the set of all finite $b$-ary strings, serving as the \emph{domain} of the structure and representing the nodes of the complete infinite $b$-ary tree;
    \item $\pr{r}$ and $\succ_i^{\mathfrak{T}}$ are defined exactly as in the case of $S3S$.
\end{itemize}

\paragraph*{Language and Syntax}
The language of $S3S$ (resp. $SbS$) is a monadic second-order logic with vocabulary $\Sig = \{\pr{c}_1, \pr{c}_2, \ldots, \pr{c}_n\} \cup \{\succ_0, \succ_1, \succ_2\}$ (resp. $\Sig = \{\pr{c}_1, \pr{c}_2, \ldots, \pr{c}_n\} \cup \{\succ_0, \succ_1,\dots, \succ_{b-1}\}$). It includes:
\begin{itemize}
\item A set of \emph{constants} denoted by $\pr{c}_1, \ldots, \pr{c}_n$.
\item \emph{Individual variables} (e.g., $x_1, x_2, x_3$) that range over elements of the domain $T$.
\item \emph{Set variables} (e.g., $X_1, X_2, X_3$) that range over subsets of $T$.
\item Logical connectives: $\wedge, \vee, \neg, \implies, \iff$.
\item Quantifiers:
\begin{itemize}
\item First-order quantifiers ($\forall x$, $\exists x$) over individual variables.
\item Second-order quantifiers ($\forall X$, $\exists X$) over set variables.
\end{itemize}
\item Atomic formulas:
\begin{itemize}
\item $x \succ_0 y$, $x \succ_1 y$ and $x \succ_2 y$ (resp. $x \succ_0 y$, $x \succ_1 y,\dots,$ and $x \succ_{b-1} y$), representing the left, middle and right successor relations.
\item $x \in X$, indicating membership of an element $x$ in a subset $X \subseteq T$.
\item $x = y$, expressing equality between individual elements.
\end{itemize}
\end{itemize}

In the following, we will also use the term ``infinite words'' to denote an infinite path of the ternary tree, starting from the root $\pr{r}$.

\section{$MSO$ inexpressibility of the hashing property}\label{sec:3}
In this section we first provide a simple proof that the property that three finite words are trifferent is not expressible in MSO logic, see Subsection \ref{subsec:33}. Then we prove, more in general, that any hashing property is not expressible in MSO logic for both finite and infinite words, see Subsections \ref{subsec:EF} and \ref{sec:32}.

\subsection{$MSO$ inexpressibility of the trifference: finite words}\label{subsec:33}
In this paragraph, we will prove, by direct proof, the nonexistence of formulas $f(X,Y,Z)$  satisfied whenever $X,Y,Z$ are finite, trifferent words.

First of all, we note that we can identify the set of finite ternary strings $T$ with the subset of binary strings $T'=\{00,01,10\}^*\subseteq B=\{0,1\}^*$ and the binary relations on $T$, $\succ_0^{\mathfrak{T}}, \succ_1^{\mathfrak{T}}, \succ_2^{\mathfrak{T}}$ with the corresponding binary relations on $T'$, $\succ_0^{\mathfrak{T'}}, \succ_1^{\mathfrak{T'}}, \succ_2^{\mathfrak{T'}}$. Here we have essentially substitute the $0$s with $00$s, the $1$s with $01$s and the $2$s with $10$s. Since the set $T'$ and the relations $\succ_0^{\mathfrak{T'}}, \succ_1^{\mathfrak{T'}}, \succ_2^{\mathfrak{T'}}$ are $S2S$ definable, if we assume, by contradiction, that the trifference relation would be $S3S$ definable, then we would obtain the existence of a formula $f'(\cdot,\cdot,\cdot)$ such that $f'(X,Y,Z)$ is true if and only if the binary finite words $X,Y,Z$ are in $T'$ and there exists a coordinate $i$, where $i$ is an odd positive integer, such that $\{X_iX_{i+1},Y_iY_{i+1},Z_iZ_{i+1}\}=\{00,01,10\}$.

We now consider the relation $R(X,Y,Z)$ on binary words that is satisfied whenever $f'(X,Y,Z)$ is true.
\begin{lem}\label{lem:S2S}
The relation $R$ is not S2S definable.
\end{lem}
\proof
Suppose, by contradiction, that the relation $R$ is $S2S$ definable. Then, according to Theorem $1(iv)$ of \cite{LauSau}, we have that, being $R$ a relation among finite words, it is a finite union of special relations. More precisely, we have
$$R=\left(\bigcup_{i\in I_1} R_i^1\right)\cup \left(\bigcup_{i\in I_2} R_i^2\right)\cup \left(\bigcup_{i\in I_3} R_i^3\right)$$
where $I_1,I_2$ and $I_3$ are finite sets and each special relation $R_i^1$ is of the form
$$R_i^1=\{(a_1a_2,a_1a_3a_4,a_1a_3a_5): a_1\in A_i^1, a_2\in A_i^2, a_3\in A_i^3, a_4\in A_i^4, a_5\in A_i^5\},$$
each $R_i^2$ is of the form
$$R_i^2=\{(a_1a_3a_4,a_1a_2,a_1a_3a_5): a_1\in A_i^1, a_2\in A_i^2, a_3\in A_i^3, a_4\in A_i^4, a_5\in A_i^5\},$$
and each $R_i^3$ is of the form
$$R_i^3=\{a_1a_3a_4,a_1a_3a_5,a_1a_2): a_1\in A_i^1, a_2\in A_i^2, a_3\in A_i^3, a_4\in A_i^4, a_5\in A_i^5\}$$
where the sets $A_i^1,\dots, A_i^5$ are regular sets. The definition of regular set is not important for this proof, but the reader can find it in \cite{LauSau}.

Given $\ell>n\geq 0$, we consider the words $X^{n,\ell},Y^{n,\ell},Z^{n,\ell}$ defined as
$$X^{n,\ell}:=(00)^{n}(01)(00)^{\ell-n-1},$$
$$Y^{n,\ell}:=(10)^{n}(00)(00)^{\ell-n-1},$$
$$Z^{n,\ell}:=(10)^{n}(10)(00)^{\ell-n-1}.$$
We have that $R(X^{n,\ell},Y^{n,\ell},Z^{n,\ell})$ is satisfied since, set $i=2n+1$, $i$ is an odd positive integer and
$$\{X^{n,\ell}_iX^{n,\ell}_{i+1},Y^{n,\ell}_iY^{n,\ell}_{i+1},Z^{n,\ell}_iZ^{n,\ell}_{i+1}\}=\{00,01,10\}.$$
This means that any triple $X^{n,\ell},Y^{n,\ell},Z^{n,\ell}$ satisfies at least one special relation of type $R_i^j$ where $j\in \{1,2,3\}$.
Therefore, fixed $\ell>|I_1|+|I_2|+|I_3|$, for the pigeonhole principle, there exists two triples $X^{n,\ell},Y^{n,\ell},Z^{n,\ell}$ and $X^{n',\ell},Y^{n',\ell},Z^{n',\ell}$ that satisfy the same relation $R_i^j$. We split the proof into two cases.

Case $1$: $X^{n,\ell},Y^{n,\ell},Z^{n,\ell}$ and $X^{n',\ell},Y^{n',\ell},Z^{n',\ell}$ both satisfy $R_i^1$. This means that
$$(X^{n,\ell},Y^{n,\ell},Z^{n,l\ell})=(a_1a_2,a_1a_3a_4,a_1a_3a_5)$$
and
$$(X^{n',\ell},Y^{n',\ell},Z^{n',\ell})=(a'_1a'_2,a'_1a'_3a'_4,a'_1a'_3a'_5).$$
Since $X^{n,\ell}$ and $Y^{n,\ell}$ begin differently, we have that $a_1$ is the emptyword, denoted here by $\underline{0}$, and
$$(X^{n,\ell},Y^{n,\ell},Z^{n,\ell})=(a_2,a_3a_4,a_3a_5).$$
Similarly, we also have that $a'_1=\underline{0}$.
Thus $X^{n,\ell}=a_2$ and $X^{n',\ell}=a'_2$ for some $a_2,a'_2$ in the same regular set $A_i^2$. But this would also imply that
$$(X^{n',\ell},Y^{n,\ell},Z^{n,\ell})=(a'_2,a_3a_4,a_3a_5)$$
satisfies the relation $R_i^1$ which is a contradiction since $(X^{n',\ell},Y^{n,\ell},Z^{n,\ell})$ does not satisfy $R$.

Case $2$: $X^{n,\ell},Y^{n,\ell},Z^{n,\ell}$ and $X^{n',\ell},Y^{n',\ell},Z^{n',\ell}$ both satisfy $R_i^2$ (or, simmetrically, $R_i^3$). This means that
$$(X^{n,\ell},Y^{n,\ell},Z^{n,\ell})=(a_1a_3a_4,a_1a_2,a_1a_3a_5)$$
and
$$(X^{n',\ell},Y^{n',\ell},Z^{n',\ell})=(a'_1a'_3a'_4,a'_1a'_2,a'_1a'_3a'_5).$$
Since $X^{n,\ell}$ and $Y^{n,\ell}$ begin differently, we have that $a_1=\underline{0}$ and
$$(X^{n,\ell},Y^{n,\ell},Z^{n,\ell})=(a_3a_4,a_2,a_3a_5).$$
$X^{n,\ell}$ and $Z^{n,\ell}$ begin differently, we have that $a_3=\underline{0}$ and
$$(X^{n,\ell},Y^{n,\ell},Z^{n,\ell})=(a_4,a_2,a_5).$$
Similarly, we also have that $a'_1=a'_3=\underline{0}$.
Thus $X^{n,\ell}=a_4$ and $X^{n',\ell}=a'_4$ for some $a_4,a'_4$ in the same regular set $A_i^4$. But this would also imply that
$$(X^{n',\ell},Y^{n,\ell},Z^{n,\ell})=(a'_4,a_2,a_5)$$
satisfies the relation $R_i^2$ which is a contradiction since $(X^{n',\ell},Y^{n,\ell},Z^{n,\ell})$ does not satisfy $R$.
\endproof
From Lemma \ref{lem:S2S}, it follows immediately that:
\begin{thm}
There is no $S3S$ formula $f(\cdot,\cdot,\cdot)$ such that $f(X,Y,Z)$ is true if and only if the finite words $X,Y,Z$ are trifferent.
\end{thm}

\subsection{Ehrenfeucht–Fraïssé game}\label{subsec:EF}
In this subsection we introduce the Ehrenfeucht–Fraïssé (EF, for short), a combinatorial game used to explore the expressive power of logical languages by comparing structures, which will play a crucial role in the next subsection. The EF game allows us to formally analyze which properties can or cannot be distinguished within a given logic, such as first-order  or monadic second-order  logic. Here, we define and analyze the EF game specifically for MSO $S3S$ logic.

The aim of the EF game on MSO $S3S$ logic is to establish whether these structures satisfy the same logical formulas up to a certain quantifier depth within this logic. We will denote this quantifier depth by $\rho$.

Given two tree structures $\mathfrak{A}$ and $\mathfrak{D}$ defined on the same domain $T$, two $\ell$-tuples $\underline{a} = (a_1, a_2, \ldots, a_{\ell}) \in T^{\ell}$ and $\underline{d} = (d_1, d_2, \ldots, d_{\ell})\in T^{\ell}$, and, finally, two $p$-tuples $\underline{V} = (V_1, V_2, \ldots, V_p)$ and $\underline{U} = (U_1, U_2, \ldots, U_p)$ where each $V_i, U_i \subseteq T$, we define the structure $\mathfrak{A}':=(\mathfrak{A}, \underline{a}, \underline{V})$ to be the structure with vocabulary $\Sig_t = \Sig \cup \{\pr{c}_1, \ldots, \pr{c}_{\ell}\} \cup \{\pr{R}_1, \ldots, \pr{R}_p\}$   where $\pr{c}_i$ is  a new constant associated to a domain elements $a_i \in T$, i.e., $\pr{c}_i^{\mathfrak{A}'} = a_i$ for every $i=1,\ldots,\ell$, and $\pr{R}_i$ is a new relation of arity $1$ associated to a set in the domain $V_i \subseteq T$, i.e., $\pr{R}_i^{\mathfrak{A}'} = V_i$ for every $i=1,\ldots,p$. Analogously, we define the structure $\mathfrak{D}':=(\mathfrak{D}, \underline{d}, \underline{U})$.

\textbf{Rules of the game.}
In the EF game, the game proceeds over $\rho$ rounds and there are two players, Alice and Bob, who play on the structures $\mathfrak{A}'$ and $\mathfrak{D}'$ with two different types of moves:
\begin{itemize}
\item \textbf{Point move.} Alice chooses a structure, $\mathfrak{A}'$ or $\mathfrak{D}'$, and an element that belongs to that structure. Bob responds with an element in the other structure.
\item \textbf{Set move.} Alice chooses a structure, $\mathfrak{A}'$ or $\mathfrak{D}'$, and a subset that belongs to that structure. Bob responds with a subset
of the other structure.
\end{itemize}

\textbf{Winning condition.}
Suppose that in a $\rho$-round game we have, for an integer $m \in \{0, \ldots, \rho\}$, the following point moves and set moves:
\begin{itemize}
\item Point moves $(x_1, \ldots, x_m) \in T^{m}$ and $(y_1, \ldots, y_m) \in T^{m}$.
\item Set moves $(X_1, \ldots, X_{\rho-m})$ and $(Y_1, \ldots, Y_{\rho-m})$ where $X_i, Y_i \subseteq T$ for each $i=1,\ldots,\rho-m$.
\end{itemize}
Bob wins the $\rho$-round game if the function
$$(\underline{a}, x_1, \ldots, x_m) \mapsto (\underline{d}, y_1, \ldots, y_m)$$ is a partial isomorphism between the two structures $(\mathfrak{A}', X_1, \ldots, X_{\rho-m})$ and $(\mathfrak{D}', Y_1, \ldots, Y_{\rho-m})$.
Alice or Bob has a winning strategy if one can guarantee that the other player will lose regardless of the move played by this player.

\begin{thm}[\cite{Libkin}]
Given two structures $\mathfrak{A}$ and $\mathfrak{D}$, for every $\rho \geq 0$, the following conditions are equivalent:
\begin{itemize}
\item Bob has a winning strategy for $\rho$-round MSO games on $\mathfrak{A}$ and $\mathfrak{D}$.
\item $\mathfrak{A} \equiv_\rho \mathfrak{D}$, i.e., $\mathfrak{A}$ and $\mathfrak{D}$ agree on MSO sentences with a number of quantifiers that does not exceed $\rho$.
\end{itemize}
\end{thm}

If a property $Q$ is true for the structure $\mathfrak{A}$, false for the structure $\mathfrak{D}$, but these two structures can be shown equivalent by providing a winning strategy for Bob for any $\rho \in \mathbb{N}$, then this shows that the property $Q$
is not expressible in MSO logic.

\subsection{$MSO$ inexpressibility of the hashing properties}\label{sec:32}
In this subsection, we consider the expressibility of the hashing properties in MSO logic. For notation convenience, we will explicitly address the case of the trifferent property on infinite words. On the other hand, with the same proof one can show that the $(b,k)$-hashing property, for finite and infinite words, is inexpressible by a single formula in MSO logic.

First, we introduce some notation. We define $T_0$ to be the subset of $T$ of the nodes which begin with $0$, $T_1$ to be the subset of the nodes which begin with $1$ and, finally, $T_2$ to be the subset of the nodes which begin with $2$. In constructing such sets we are assuming that the root $\pr{r}$ belongs to $T_0$, $T_1$ and $T_2$.

Given a structure $(\mathfrak{T}, x_1,\dots,x_{\ell},X_0,\dots,X_{h-\ell})$ we define its restriction to $T_j$, for each $j=0,1,2$, as
$$(\mathfrak{T}, x_1,\dots,x_{\ell},X_0,\dots,X_{h-\ell})|_{T_j}=(\mathfrak{T}|_{T_j}, x_i: x_i\in T_j,X_0\cap T_j,\dots,X_{h-\ell}\cap T_j)$$
where $$\mathfrak{T}|_{T_j}=\left(T_j, (\succ_0^{\mathfrak{T}})|_{T_j}, (\succ_1^{\mathfrak{T}})|_{T_j}, (\succ_2^{\mathfrak{T}})|_{T_j}\right).$$

\begin{lem}\label{xAndy}
Given $M\in \mathbb{N}$, in the ternary $\mathfrak{T}$ tree, there exist two different infinite words of the form $X_0=0\dots0100\dots$ and $Y_0=0\dots0100\dots$ (where the number of zeros before the one differs in the two cases) and, for any formula $p(\cdot)$  of length at most $M$, it results $p(X_0)\iff p(Y_0)$.
\end{lem}
\proof
Denoted by $\mathcal{F}_M$ the family of the formulas of length at most $M$ with exactly one free variable (of set type), we have that $|\mathcal{F}_M|$ is finite. On the other hand, $T_0$ contains infinitely many infinite words of the form $0\dots0100\dots$. It follows that there are infinitely many infinite words of that form that coincide in any formulas of $\mathcal{F}_M$. Hence we can also find two infinite words $X_0,Y_0 \subseteq T_0$ such that $X_0=0\dots0100\dots$ and $Y_0=0\dots0100\dots$ and the number of zeros before the one differs in the two cases and such that for any formula $p \in \mathcal{F}_M$, it holds $p(X_0)\iff p(Y_0)$.
\endproof

\begin{rem}
We have already noted that the set of all the formulas of length at most $M$ with exactly one free set variable (resp. individual variable) is finite. The same can also be said for the family $\mathcal{F}_{M,i,j}$ of the formulas of length at most $M$ with exactly $i$ free variables of individual type and $j$ free variables of set type.
\end{rem}

\begin{prop}\label{AutoT0}
There exist two infinite words $X_0,Y_0\subseteq T_0$ of the form $X_0=0\dots0100\dots$ and $Y_0=0\dots0100\dots$ (with a different number of zeros before the one) such that $(\mathfrak{T},X_0)|_{T_0}\equiv_\rho (\mathfrak{T},Y_0)|_{T_0}$.
\end{prop}
\proof
Let us consider a sequence $M_0,M_1,\dots,M_{\rho}$ of natural numbers such that $M_0=1$ and, for any $h$, $1 \leq h\leq \rho$, take
$$M_h\geq \sum_{i,j\leq \rho}|\mathcal{F}_{M_{h-1},i,j}|(M_{h-1}+2)+2.$$
Note that such a sequence exists since $|\mathcal{F}_{M_{h-1},i,j}|$ is finite.

According to Lemma \ref{xAndy}, we can choose $X_0=0\dots0100\dots$ and $Y_0=0\dots0100\dots$, where the number of zeros before the one differs in the two cases, be two infinite words of $T_0$ such that for any formula $p(\cdot)$ of length at most $M_\rho$ we have $p(X_0)\iff p(Y_0)$.

We prove that $(\mathfrak{T},X_0)|_{T_0}\equiv_\rho (\mathfrak{T},Y_0)|_{T_0}$ inductively, using the EF game. More precisely, we prove that, if at the $h$-th step Alice chooses the point $x_{\ell}$ (resp. the set $X_{h-\ell+1}$) and considering the structure $(x_1,\dots, x_\ell,X_0,\dots,$ $X_{h-\ell})$, Bob can choose $y_{\ell}$ (resp. the set $Y_{h-\ell+1}$) so that any formula of $\mathcal{F}_{M_{\rho-h},\ell,h-\ell+1}$ is realized on $(y_1,\dots, y_\ell,Y_0,\dots,Y_{h-\ell})$ if and only if it holds on $(x_1,\dots, x_\ell,X_0,\dots,X_{h-\ell})$. Equivalently we can say that $(y_1,\dots, y_\ell,Y_0,\dots,$ $Y_{h-\ell})$ is equivalent to $(x_1,\dots, x_\ell,X_0,\dots,X_{h-\ell})$ for any formula of length at most $M_{\rho-h}$.

\textbf{Base step:}

Let us suppose that Alice chooses a node $x_1$.
We denote by $\mathcal{F}^t_{M_{\rho-1},1,1}$ the set of all the formulas of $\mathcal{F}_{M_{\rho-1},1,1}$ that are true on $(x_1, X_0)$ and we set
$$\mathcal{F}^f_{M_{\rho-1},1,1}:=\mathcal{F}_{M_{\rho-1},1,1}\setminus \mathcal{F}^t_{M_{\rho-1},1,1}.$$

We note that
$$F_\rho(\cdot):=\exists z: \cup_{f\in \mathcal{F}^t_{M_{\rho-1},1,1}} f(z, \cdot)\cup_{g\in \mathcal{F}^f_{M_{\rho-1},1,1}} \neg g(z, \cdot)$$

has length at most $M_{\rho}$. Also, this formula is true on $X_0$ and so this means that it is true also on $Y_0$, but then
$$\exists z:\cup_{f\in \mathcal{F}^t_{M_{\rho-1},1,1}} f(z, Y_0)\cup_{g\in \mathcal{F}^f_{M_{\rho-1},1,1}} \neg g(z, Y_0).$$

Hence Bob can choose $y_1$ so that for any formula $p\in \mathcal{F}_{M_{\rho-1},1,1}$, it holds: $p(x_1, X_0)\iff p(y_1, Y_0).$

Similarly, if Alice chooses a set $X_1$, Bob can find a set $Y_1$ so that for any formula $p\in \mathcal{F}_{M_{\rho-1},0,2}$, we have $p(X_0,X_1)\iff p(Y_0,Y_1)$.

\textbf{Inductive Step:}

Let us suppose that, at the $h$-th step, with $h \leq \rho$, Alice chooses a node $x_\ell$.
We denote by $\mathcal{F}^t_{M_{\rho-h},\ell,h-\ell+1}$ the set of all the formulas of $\mathcal{F}_{M_{\rho-h},\ell,h-\ell+1}$ that are true on $(x_1,\dots, x_\ell,X_0,\dots X_{h-\ell})$ and we set
$$\mathcal{F}^f_{M_{\rho-h},\ell,h-\ell+1}:=\mathcal{F}_{M_{\rho-h},\ell,h-\ell+1}\setminus \mathcal{F}^t_{M_{\rho-h},\ell,h-\ell+1}.$$

Then, proceeding as before, Bob chooses $y_\ell$ so that for any formula $p\in \mathcal{F}_{M_{\rho-h},\ell,h-\ell+1}$, we have:
$$p(x_1,\dots, x_\ell,X_0,\dots X_{h-\ell})\iff p(y_1,\dots, y_\ell,Y_0,\dots Y_{h-\ell}).$$ 
Note that such a $y_\ell$ exists because we can consider the formula
\begin{multline*}
F_{\rho-(h-1)}(\cdot):= \\ \exists z:\cup_{f\in \mathcal{F}^t_{M_{\rho-h},\ell,h-\ell+1}} f(\cdots,z,\cdots)\cup_{g\in \mathcal{F}^f_{M_{\rho-h},\ell,h-\ell+1}} \neg g(\cdots,z,\cdots)
\end{multline*}
which is true on $(x_1,\dots,x_{\ell-1},X_0,\dots,X_{h-\ell})$ and so, since its length is smaller than $M_{\rho-(h-1)}$, it must be true also on $$(y_1,\dots,y_{\ell-1},Y_0,\dots,Y_{h-\ell}).$$

Similarly, if Alice chooses a set $X_{h-\ell+1}$, Bob can find a set $Y_{h-\ell+1}$ so that for any formula $p\in \mathcal{F}_{M_{\rho-h},\ell-1,h-\ell+2}$, we have:
$$p(x_1,\dots, x_{\ell-1},X_0,\dots, X_{h-\ell+1})\iff p(y_1,\dots, y_{\ell-1},Y_0,\dots, Y_{h-\ell+1}).$$
\endproof

\begin{prop}\label{Incollamento}
Assume $\rho \geq h\geq \ell \geq 1$ and $$(\mathfrak{T},x_1,\dots,x_{\ell-1},X_0,\dots, X_{h-\ell})|_{T_j}\equiv_\rho (\mathfrak{T},y_1,\dots,y_{\ell-1},Y_0,\dots, Y_{h-\ell})|_{T_j},$$
for $j = 0,1,2$.
This implies that
$$(\mathfrak{T},x_1,\dots,x_{\ell-1},X_0,\dots, X_{h-\ell})\equiv_\rho (\mathfrak{T},y_1,\dots,y_{\ell-1},Y_0,\dots, Y_{h-\ell}).$$
\end{prop}
\proof
Note that if $\ell=1$, no $x_i$ and $y_i$ appears in the structure.
We prove the thesis
using the EF game inductively on $\rho$.

\textbf{Base step}. \\
Let us suppose, by contradiction, that
$$(\mathfrak{T},x_1,\dots,x_{\ell-1},X_0,\dots, X_{h-\ell})\not\equiv_0 (\mathfrak{T},y_1,\dots,y_{\ell-1},Y_0,\dots, Y_{h-\ell}).$$
This means that either $x_j\succ_0 x_i$ (resp. $\succ_1$, $\succ_2$) but $y_j\not\succ_0 y_i$ (resp. $\succ_1$, $\succ_2$) or $x_j \in X_i$ but $y_j\not\in Y_i$.
In the first case, since $x_j\succ_0 x_i$, we may assume that $x_j,x_i \in T_0$. Hence, since
$$(\mathfrak{T},x_1,\dots,x_{\ell-1},X_0,\dots, X_{h-\ell})|_{T_0}\equiv_0 (\mathfrak{T},y_1,\dots,y_{\ell-1},Y_0,\dots, Y_{h-\ell})|_{T_0}$$
we also have $y_j,y_i\in T_0$ and $y_j\succ_0 y_i$.

We proceed similarly in the second case. Here we may assume $x_j\in T_0$ and $x_j\in X_i\cap T_0$ and it would follow that also $y_j\in T_0$, and $y_j\in Y_i\cap T_0$. Hence we must have that
$$(\mathfrak{T},x_1,\dots,x_{\ell-1},X_0,\dots, X_{h-\ell})\equiv_0 (\mathfrak{T},y_1,\dots,y_{\ell-1},Y_0,\dots, Y_{h-\ell}).$$

\textbf{Inductive step.} \\
Let us suppose that Alice chooses a node $x_{\ell+1}\in T_0$.
Since $$(\mathfrak{T},x_1,\dots,x_{\ell},X_0,\dots, X_{h-\ell})|_{T_0}\equiv_\rho (\mathfrak{T},y_1,\dots,y_{\ell},Y_0,\dots, Y_{h-\ell})|_{T_0},$$ Bob can choose $y_{\ell+1}\in T_0$ so that
$$(\mathfrak{T},x_1,\dots,x_{\ell},X_0,\dots, X_{h-\ell})|_{T_0}\equiv_{\rho-1} (\mathfrak{T},y_1,\dots,y_{\ell},Y_0,\dots, Y_{h-\ell})|_{T_0}.$$
Since, for $j=1,2$,
$$(\mathfrak{T},x_1,\dots,x_{\ell},X_0,\dots, X_{h-\ell})|_{T_j}=(\mathfrak{T},x_1,\dots,x_{\ell-1},X_0,\dots, X_{h-\ell})|_{T_j}$$
and
$$ (\mathfrak{T},y_1,\dots,y_{\ell},Y_0,\dots, Y_{h-\ell})|_{T_j}=(\mathfrak{T},y_1,\dots,y_{\ell-1},Y_0,\dots, Y_{h-\ell})|_{T_j}.$$
We also have that, for $j=1,2$,
$$(\mathfrak{T},x_1,\dots,x_{\ell},X_0,\dots, X_{h-\ell})|_{T_j}\equiv_{\rho-1} (\mathfrak{T},y_1,\dots,y_{\ell},Y_0,\dots, Y_{h-\ell})|_{T_j}.$$
Due to the inductive hypothesis, we obtain that
$$(\mathfrak{T},x_1,\dots,x_{\ell},X_0,\dots, X_{h-\ell})\equiv_{\rho-1} (\mathfrak{T},y_1,\dots,y_{\ell},Y_0,\dots, Y_{h-\ell}).$$

But this means that Bob has a winning strategy for the $\rho-1$ moves EF game on
$$(\mathfrak{T},x_1,\dots,x_{\ell},X_0,\dots, X_{h-\ell})\mbox{ and } (\mathfrak{T},y_1,\dots,y_{\ell},Y_0,\dots, Y_{h-\ell}).$$

If Alice chooses a set $X_{h-\ell+1}$, proceeding similarly, Bob can choose sets $Y^0_{h-\ell+1}\subseteq T_0$, $Y^1_{h-\ell+1}\subseteq T_1$ and $Y^2_{h-\ell+1}\subseteq T_2$ such that 
$$(\mathfrak{T},x_1,\dots,x_{\ell-1},X_0,\dots, X_{h-\ell+1})|_{T_j}\equiv_{\rho-1} (\mathfrak{T},y_1,\dots,y_{\ell-1},Y_0,\dots, Y^j_{h-\ell+1})|_{T_j}$$
for $j=0,1,2$.
Setting 
$$Y_{h-\ell+1}:=Y^0_{h-\ell+1}\cup Y^1_{h-\ell+1} \cup Y^2_{h-\ell+1} $$
we obtain that 
$$(\mathfrak{T},x_1,\dots,x_{\ell-1},X_0,\dots, X_{h-\ell+1})\equiv_{\rho-1} (\mathfrak{T},y_1,\dots,y_{\ell-1},Y_0,\dots, Y_{h-\ell+1}).$$
It follows that Bob has a winning strategy for the $\rho-1$ moves EF game on
$$(\mathfrak{T},x_1,\dots,x_{\ell-1},X_0,\dots, X_{h-\ell}, X_{h-\ell+1})$$ and  $$(\mathfrak{T},y_1,\dots,y_{\ell-1},Y_0,\dots, Y_{h-\ell}, Y_{h-\ell+1}).$$

Summing up these two cases, we obtain that Bob has a winning strategy for the $\rho$ moves EF game on
$$(\mathfrak{T},x_1,\dots,x_{\ell-1},X_0,\dots, X_{h-\ell})$$ and  $$(\mathfrak{T},y_1,\dots,y_{\ell-1},Y_0,\dots, Y_{h-\ell}).$$
\endproof

\begin{thm}
There is no $S3S$ formula $f(\cdot,\cdot,\cdot)$ such that $f(X,Y,Z)$ is true if and only if the infinite words $X,Y,Z$ are trifferent.
\end{thm}
\proof
Let us assume, by contradiction, that such a formula $f$ exists and let $\rho$ be its rank.
We consider distinct infinite words $X_0,Y_0\subseteq T_0$ of the form $0\dots0100\dots$ such that $(\mathfrak{T},X_0)|_{T_0}\equiv_\rho (\mathfrak{T},Y_0)|_{T_0}$. Let now consider $X_1, X_2\subseteq T_1$ of the form $X_1=1\dots111\dots$ and $X_2=1\dots1211\dots$ where the number $m$ of ones before the $2$ in $X_2$ coincides with the number of zeros before the $1$ of $X_0$.
Set $Y_1=X_1$ and $Y_2=X_2$ we have that,
$$(\mathfrak{T},X_0)|_{T_0}\equiv_\rho (\mathfrak{T},Y_0)|_{T_0}$$
and also
$$(\mathfrak{T},X_1,X_2)|_{T_1}\equiv_\rho (\mathfrak{T},Y_1,Y_2)|_{T_1}.$$
Hence, according to Proposition \ref{Incollamento}, we have that
$$(\mathfrak{T},X_0,X_1,X_2)\equiv_\rho (\mathfrak{T},Y_0,Y_1,Y_2).$$
On the other hand the words $X_0,X_1,X_2$ are not trifferent and hence
$f(X_0,$ $X_1,X_2)$ is not realized. Then, we also have that $f(Y_0,Y_1,Y_2)$ should be false.
But, since in the $(m+1)$-th coordinate of $Y_0,Y_1,Y_2$ are, respectively, $0, 1, 2$, these words are trifferent, contradicting the hypothesis that $f(X,Y,Z)$ is true if and only if $X,Y,Z$ are trifferent.
\endproof

With the same proof, one can check that the following result holds.
\begin{thm}
Let $3\leq k\leq b$. Then, there is no $SbS$ formula $f(\cdot,\cdot,\ldots,\cdot)$ such that $f(X_0,X_1,\ldots,$ $ X_{k-1})$ is true if and only if the infinite words $(X_0,X_1,\ldots$ $, X_{k-1})$ satisfy the $(b,k)$-hashing property.

Analogously, there is no $SbS$ formula $f(\cdot,\cdot,\ldots,\cdot)$ such that $f(x_1,x_2,\ldots, x_k)$ is true if and only if the finite words $(x_1,x_2,\ldots,x_k)$ satisfy the $(b,k)$-hashing property.
\end{thm}

\section{A computational approach}
In this section, we use the Z3 Satisfiability Modulo Theories solver (SMT solver, for short) \cite{Z3} developed by Microsoft Research. More precisely we write our code in Z3Py which is a Python binding for Z3, enabling intuitive modeling and solving of complex logical constraints in Python.

 We are interested in studying the existence of a trifferent code $C$ of \(m\) codewords, each of length \(n\), over a ternary alphabet \(\{0, 1, 2\}\). In particular, we run the following program in the case $m=11$ and $n=5$ obtaining the non-existence of such code (recovering a result of \cite{DFGP}).

Here we provide an implementation of the trifference property in Z3Py.
\small
\begin{lstlisting}[language=Python, label=lst:trifference]
from z3 import *
import itertools

def findsubsets(S, k):
    return set(itertools.combinations(S, k))

def getTrifferenceSolver(m, n):
    C = [[ Int("x_%s_%s" % (i+1, j+1)) for j in range(n) ]
          for i in range(m) ]
    triplet_indexes = findsubsets(range(m), 3)
    trifference_property = [Or([ 
        And(C[t[0]][i] != C[t[1]][i], 
            C[t[0]][i] != C[t[2]][i], 
            C[t[1]][i] != C[t[2]][i]) 
        for i in range(n) ]) for t in triplet_indexes]
    
    vocabulary = [[ And(0 <= C[i][j], C[i][j] <= 2) 
                   for j in range(n)] 
                   for i in range(m)]
    vocabulary = sum(vocabulary, [])
    s = Solver()
    s.add(trifference_property)
    s.add(vocabulary)

    return [s, C]
\end{lstlisting}
\normalsize
The function \texttt{getTrifferenceSolver(m, n)} builds a matrix \(C\) where each row corresponds to a codeword. The constraints added to the SMT solver are as follows:

\begin{itemize}
    \item \textbf{Trifference Property:} For every triplet of rows, at least one column must have values that are pairwise distinct.
    \item \textbf{Vocabulary Constraint:} Each entry in the matrix is constrained to the ternary alphabet \(\{0, 1, 2\}\).
\end{itemize}

The solver could be then invoked to check the satisfiability of these constraints. If a solution exists, the values of \(C\) can be extracted, otherwise, a message indicating unsatisfiability can be displayed.

\section*{Acknowledgements}
This work was supported by the European Union under the Italian National Recovery and Resilience Plan (NRRP) of NextGenerationEU, Partnership on “Telecommunications of the Future,” Program “RESTART” under Grant PE00000001, “Netwin” Project (CUP E83C22004640001).\\


\begin{thebibliography}{99}
\bibitem{B} Bhandari S. and Khetan A., \emph{Improved Upper Bound for the Size of a Trifferent Code}, Combinatorica 45, 2 (2025).
\bibitem{bis-dha-gij-2024} Bishnoi, A. and D'haeseleer, J. and Gijswijt, D. and Potukuchi, A., \emph{Blocking sets, minimal codes and trifferent codes}, Journal of the London Mathematical Society, Vol. 109, No. 6, (2024).
\bibitem{CD} Costa S. and Dalai M., \emph{A gap in the slice rank of k-tensors}, Journal of Combinatorial Theory, Series A 177 (2021), 105335.
\bibitem{C1} Courcelle B., \emph{The monadic second-order logic of graphs. I. Recognizable sets of finite graphs}, Information and Computation, 85 (1990), 12--75.
\bibitem{C2} Courcelle B., \emph{The monadic second-order logic of graphs, II: Infinite graphs of bounded
width}, Mathematical Systems Theory, 21 (1988), 187--221.
\bibitem{DFGP} Della Fiore S., Gnutti A. and Polak S., \emph{The maximum cardinality of trifferent codes with lengths 5 and 6}, Examples and Counterexamples 2 (2022), 100051.
\bibitem{Z3} De Moura L. and Bjørner N., \emph{Z3: An efficient SMT solver}, International conference on Tools and Algorithms for the Construction and Analysis of Systems, Springer Berlin Heidelberg, 2008.
\bibitem{Mona}  Henriksen J. G., Henriksen, Jensen J., Jørgensen M., Klarlund N., Paige R., Rauhe T. and Sandholm A., \emph{Mona: Monadic second-order logic in practice}, Tools and Algorithms for the Construction and Analysis of Systems: First International Workshop, TACAS'95 Aarhus, Denmark, May 19–20, 1995 Selected Papers 1. Springer Berlin Heidelberg, 1995.
\bibitem{Felix} Klaedtke F. and Rueß H., \emph{Monadic second-order logics with cardinalities}, 30th International Colloquium (ICALP 2003), 2719 (2003), Springer, 681--696.
\bibitem{KoS} K\H{o}rner J. and Gábor S., \emph{Trifference}, Studia Scientiarum Mathematicarum Hungarica 30 (1995), 95--103.
\bibitem{KM} K\H{o}rner J. and Marton K., \emph{New bounds for perfect hashing via information theory}, European Journal of Combinatorics 9.6 (1988), 523--530.
\bibitem{K1} Kuncak V., Nguyen H. H. and Rinard M., \emph{An Algorithm for Deciding BAPA:
Boolean Algebra with Presburger Arithmetic}, Automated Deduction
(CADE 2005), 3632 (2005), 260--277, Springer, 2005.
\bibitem{K2} Kuncak V., Nguyen H. H. and Rinard M., \emph{Deciding Boolean Algebra with
Presburger Arithmetic}, Journal of Automated Reasoning, 36 (2006), 213--239.
\bibitem{KS} Kurz S., \emph{Trifferent codes with small lengths}, Examples and Counterexamples 5 (2024), 100139.
\bibitem{LauSau} Lauchli H. and Savioz C. \emph{Monadic Second Order Definable Relations on the Binary Tree}, The Journal of Symbolic Logic, Vol. 52, No. 1 (1987), 219-226.
\bibitem{Libkin} Libkin L., \emph{Elements of finite model theory}, Vol. 41. Heidelberg: Springer, 2004.
\bibitem{Rudolph} Herrmann L., Peth V. and Rudolph S., \emph{Decidable (ac) counting with Parikh and Muller: adding Presburger arithmetic to monadic second-order logic over tree-interpretable structures}, 32nd EACSL Annual Conference on Computer Science Logic, 2024.
\end{thebibliography}
\end{document}